\documentclass[11pt]{amsart}

\usepackage{amssymb}
\usepackage{amsmath}
\usepackage{amsthm}
\usepackage{stmaryrd}
\usepackage{amsfonts,mathrsfs}
\usepackage{fullpage}
\usepackage{color}                    
\usepackage{hyperref}                 

\title{Higher cohomology of the pluricanonical bundle is not deformation invariant}
\author{Ning Hao, Li Li}


\begin{document}
\maketitle

\theoremstyle{plain}
\newtheorem{thm}{\sc Theorem}
\newtheorem{lem}{\sc Lemma}[section]
\newtheorem{d-thm}[lem]{\sc Theorem}
\newtheorem{prop}[lem]{\sc Proposition}
\newtheorem{cor}[lem]{\sc Corollary}

\theoremstyle{definition}
\newtheorem{conj}[lem]{\sc Conjecture}
\newtheorem{defn}[lem]{\sc Definition}
\newtheorem{ques}[lem]{\sc Question}

\theoremstyle{definition}
\newtheorem{ex}{\sc Example}

\theoremstyle{remark}
\newtheorem*{rmk}{\sc Remark}
\newtheorem*{rmks}{\sc Remarks}
\newtheorem*{ack}{\sc Acknowledgment}

\begin{abstract}
In this paper, we compute the dimensions of the $1$st and $2$nd
cohomology groups of all the pluricanonical bundles for Hirzebruch
surfaces, and the dimensions of the $1$st and $2$nd cohomology
groups of the second pluricanonical bundles for a blow-up of a
projective plane along finite many distinct points. Both of them
give examples showing that the higher cohomology of the
pluricanonical bundle is not deformation invariant.
\end{abstract}

{\bf \tableofcontents}

\section{Introduction}

Let $\Delta=\{t \in\mathbb{C}: |t|<1\}$ be the unit disc in the complex plane.  Assume $Y$ is a
compact complex manifold. Denote by $K_Y$ the canonical line bundle of $Y$. Denote by $h^i(Y,
L)=\dim_\mathbb{C} H^i(Y, L)$ for any line bundle $L$ on $Y$. The dimension of $H^0(Y, mK_Y)$ is
called the $m$-th plurigenus of $Y$. In \cite{Siu}, Siu proved the following remarkable theorem on
the deformation invariance of the plurigenera.

\begin{thm}\label{plurigenera}\cite{Siu}
Let $\pi: X\to \Delta$ be a holomorphic family of compact complex
projective algebraic manifolds with fiber $X_t$. Then for any
positive integer $m$, the complex dimension of $H^0(X_t,
mK_{X_t})$ is independent of $t$ for $t\in \Delta$.
\end{thm}

It is then natural to ask the following question:
\begin{ques}\label{q1}
Is the dimension of $H^q(X_t, mK_{X_t})$ independent of $t$, for
$q>0$?
\end{ques}

The main point of the present article is to give a negative answer to Question \ref{q1}.  In fact, the following two theorems are our main results.

\begin{thm}\label{Hirzebruch}
For $m\ge 0$, Let $F_m=\mathbb{P}(\mathcal{O}\oplus\mathcal{O}(m))$
be the $m$-th Hirzebruch surface. Then
\[
\dim H^0(-kK_{F_m})=\frac{1}{2}\Big{(}4k+(k-\Big{[}\frac{2k}{m}\Big{]})m+2\Big{)}\Big{(}k+\Big{[}\frac{2k}{m}\Big{]}+1\Big{)}.
\]
\end{thm}

\begin{thm}\label{blow-up}
Let $M\to\mathbb{P}^2$ be the blow-up of $\mathbb{P}^2$ along $v$ distinct points. The dimension of $H^0(-K_M)$ is $10-v$ for $v\le 4$; for $v>4$, the dimension depends on the position of the $v$ points: for each integer $i$ such that $\max(10-v,0)\le i\le6$, there exist $v$ distinct points such that the corresponding $M$ satisfies $\dim H^0(-K_M)=i$.
\end{thm}

Returning to Question \ref{q1}, Theorem \ref{Hirzebruch} shows that for each $k$ the dimensions $h^0(-kK_{F_m})$ depend on $m$ in a non-constant way.  On the other hand, it is shown in \cite[Example 2.16]{Kod} that for any positive integer $\ell \le m/2$ there is a holomorphic family $M \to \Delta$ whose fiber $M_t$ is isomorphic to $F_{m-2\ell}$ for $t \neq 0$, and whose central fiber $M_0$ is isomorphic to $F_m$.  By Serre duality, we find that
\[
h^2((k+1)K_{M_t}) = h^0(-kK_{M_t})
\]
depends on $t$ in a non-constant fashion.

Similarly, consider the situation of Theorem \ref{blow-up}.  Let $B=\{(q_1,\cdots,q_v)\in (\mathbb{P}^2)^v: q_i\neq q_j
\textrm{ if } i\neq j\}$ be the configuration space of $v$ points in
$\mathbb{P}^2$. Denote by $\Gamma_{\pi_i}\subset B\times
\mathbb{P}^2$ the graph of the morphism $\pi_i: B\to\mathbb{P}^2$
induced from the projection $(\mathbb{P}^2)^v\to \mathbb{P}^2$ to
the $i$-th factor. Let $\widetilde{M}$ be the blow-up of $B\times\mathbb{P}^2$ along $v$ disjoint submanifolds $\Gamma_{\pi_i}$ for $1\le i\le v$. The natural morphism $\widetilde{M}\to B$ is a holomorphic family with each fiber $M_t$ over $t=(q_1,\cdots,q_v)$ being isomorphic to the blow-up of $\mathbb{P}^2$ along $v$ points $q_1,\cdots,q_v$.  Theorem \ref{blow-up} says that
\[
h^2(2K_{M_t}) = h^0(-K_{M_t})
\]
is not constant for $v>5$.  Since $B$ is obviously smooth and connected, we can find a small disc
$\Delta\subset B$ such  that $h^2(2K_{M_t}))$ is not constant for $t\in \Delta$.

\begin{rmk}
So far we have shown, modulo the proofs of Theorems \ref{Hirzebruch} and \ref{blow-up}, that there are holomorphic families $M \to \Delta$ such that $h^2(kK_{M_t})$ depends on $t$.  On the other hand, we know from Theorem \ref{plurigenera} (and in fact from earlier results using the classification of surfaces) that $h^0(kK_{M_t})$ {\it is} constant in $t$.  Now, suppose $X$ is a complex projective algebraic surface. Let $K$ be the canonical line bundle of $X$. Applying the Riemann-Roch Theorem for surfaces, we have
$$
h^0(kK)-h^1(kK)+h^2(kK)=\frac{1}{2}kK(kK-K)+\chi(\mathcal{O}_X)=\frac{1}{2}k(k-1)K^2+\chi(\mathcal{O}_X).
$$
Moreover, by Noether's formula,
$$
\chi(\mathcal{O}_X)=\frac{K^2+\chi_{top}(X)}{12},
$$
and therefore
$$
h^1(kK)=h^0(kK)+h^2(kK)-(6k^2-6k+1)\chi(\mathcal{O}_X)+\frac{1}{2}k(k-1)\chi_{top}(X).
$$
Now, a deformation does not change the differentiable structure,
and $\chi_{top}(X)$ is a topological invariant. Therefore,
$\chi_{top}(X)$ is a deformation invariant.  Next,
$\chi(\mathcal{O}_X)=h^{0,0}-h^{1,0}+h^{2,0}$ is also invariant
for a family of K\"ahler manifolds, because Hodge numbers are
invariant for such a family. By Siu's result of invariance of
plurigenera, $h^0(kK)$ is also a deformation invariant.  It
follows that, since $h^2(kK)$ is not a deformation invariant,
neither is $h^1(kK)$.
\end{rmk}

\noindent{Acknowledgements.} We are grateful to Dror Varolin, who brought this problem to our
attention and provided considerable advice and encouragement throughout our investigations, to
Blaine Lawson for his valuable comments and encouragement, to Yujen Shu for useful discussions.

\section{Deformation of a Hirzebruch surface: Proof of Theorem \ref{Hirzebruch}}
Recall that the Hirzebruch surfaces $F_m$ ($m\ge 0$) is the
$\mathbb{P}^1$-bundle over $\mathbb{P}^1$ associated to the sheaf
$\mathcal{O}_{\mathbb{P}^1}\oplus\mathcal{O}_{\mathbb{P}^1}(m)$. For example, $F_0$ is isomorphic
to $\mathbb{P}^1\times \mathbb{P}^1$, and $F_1$ is the blow-up of
$\mathbb{P}^2$ at a point.

$\mathbb{P}^1$ can be covered by two affine line $U_1, U_2\cong
\mathbb{C}$. Let $z_1$ and $z_2$ be affine coordinates for $U_1$ and $U_2$ respectively. Then on $U_1$ and $U_2$ we have $z_1z_2=1$.  The sheaf $\mathcal{O}_{\mathbb{P}^1}(m)$ can be constructed by gluing two trivial bundles $U_1\times \mathbb{C}$ and $U_2\times \mathbb{C}$ together as follows: a section of $s_1$ of $U_1\times \mathbb{C}$ can be glued with a section $s_2$ of $U_2\times \mathbb{C}$ if $s_1(z_1)=z_1^ms_2(z_2)$ on $U_1\cap U_2$ where $z_2=1/z_1$ (that is to say, the transition function of the corresponding vector bundle is $z_1^m$). Therefore, $F_m$ can be obtained by gluing $U_1\times \mathbb{P}^1$ with $U_2\times \mathbb{P}^1$ where $(z_1, \zeta_1)\in U_1\times \mathbb{P}^1$ is identified with $(z_2, \zeta_2)\in U_2\times\mathbb{P}^1$ if
\[
z_1z_2=1, \textrm{ and } \zeta_1=z_2^m\zeta_2.
\]
(Here and below we adopt the standard abuse of notation, using an affine coordinate $\zeta _1$ to denote a point in the fiber ${\mathbb{P}^1}$.)

We recall the holomorphic family $\{M_t: t\in\mathbb{C}\}$ defined in \cite[Example 2.16]{Kod}: Fix a positive integer $\ell \le m/2$, define $M_t$ as
\[
M_t=U_1\times\mathbb{P}^1\cup U_2\times\mathbb{P}^1,
\]
where $(z_1,\zeta_1)\in U_1\times \mathbb{P}^1$ and $(z_2,\zeta_2)\in U_2\times \mathbb{P}^1$ are the same point in $M_t$ if
\[
z_1z_2=1 \textrm{ and } \zeta_1=z_2^m\zeta_2+tz_2^{\ell}.
\]
It is proved in \cite{Kod} that $M_t$ is isomorphic to Hirzebruch surface $F_m$ if $t=0$ and
isomorphic to $F_{m-2\ell}$ if $t\neq 0$. Thus $F_m$ is a deformation of $F_{m-2\ell}$.

Now we calculate $\dim H^0(-kK_{F_m})=\dim \Gamma  \big{(}\big{(} \bigwedge
^2T_{F_m}\big{)}^{\bigotimes k}\big{)}$ for $k>0$.  First, we assert that a section in
$\Gamma\big{(}U_1\times \mathbb{P}^1,\big{(}\bigwedge\nolimits^2 T_{U_1\times
\mathbb{P}^1}\big{)}^{\otimes k}\big{)}$ is of the form
$$v_1(z_1,\zeta_1)=\Big{(}a_{2k}(z_1)\zeta_1^{2k}+a_{2k-1}(z_1)\zeta_1^{2k-1}+...+a_0(z_1)\Big{)}\Big{(}\frac{\partial}{\partial
z_1}\wedge\frac{\partial}{\partial \zeta_1}\Big{)}^{\otimes k},$$ i.e., it the expansion contains
no term of power of $\zeta^i$ for $i>2k$. Indeed, expand $v_1(z_1,\zeta_1)$ into a power series
$\big{(}\sum_{i=0}^\infty a_i(z_1)\zeta_1^i\big{)}(\frac{\partial}{\partial
z_1}\wedge\frac{\partial}{\partial \zeta_1})^{\otimes k}$. Now cover $U_1\times \mathbb{P}^1$ with
two open subset $U_1\times V_1$ and $U_1\times V_2$ where $V_1$ and $V_2$ are isomorphic to
$\mathbb{C}$ and cover $\mathbb{P}^1$. Assume the coordinates of $V_1$ and $V_2$ are $\zeta_1$ and
$\eta_1$ respectively. Then $\zeta_1\eta_1=1$ on $V_1\cap V_2$, and therefore
$\frac{\partial}{\partial \zeta_1}=-\eta_1^2\frac{\partial}{\partial \eta_1}$. In $U_1\times V_2$,
the section $v_1(z_1,\zeta_1)$ is of the form
\[
\Big{(}\sum_{i=0}^\infty a_i(z_1)\frac{1}{\eta_1^i}\Big{)}\Big{(}\frac{\partial}{\partial
z_1}\wedge(-\eta_1^2)\frac{\partial}{\partial \eta_1}\Big{)}^{\otimes
k}=(-1)^k\Big{(}\sum_{i=0}^\infty a_i(z_1){\eta_1^{2k-i}}\Big{)}\Big{(}\frac{\partial}{\partial
z_1}\wedge\frac{\partial}{\partial \eta_1}\Big{)}^{\otimes k}.
\]
For this section to be holomorphic, all the terms $a_i(z_1)$ for $i>2k$
should vanish.

Similarly, a section in $\Gamma\big{(}U_2\times \mathbb{P}^1,\big{(}\bigwedge\nolimits^2
T_{U_2\times \mathbb{P}^1}\big{)}^{\otimes k}\big{)}$ is of the form
\[
v_2(z_2,\zeta_2)=\Big{(}b_{2k}(z_2)\zeta_2^{2k}+b_{2k-1}(z_2)\zeta_2^{2k-1}+...+b_0(z_2)\Big{)}\Big{(}\frac{\partial}{\partial
z_2}\wedge\frac{\partial}{\partial \zeta_2}\Big{)}^{\otimes k}.
\]
To glue $v_1$ and $v_2$ into a section of $\Gamma\big{(}F_m,\big{(}\bigwedge\nolimits^2
T_{F_m}\big{)}^{\otimes k}\big{)}$,  we need $v_1(z_1,\zeta_1)=v_2(z_2,\zeta_2)$ over $(U_1 \bigcap
U_2) \times \mathbb{P}^1$. Since
\[
\frac{\partial}{\partial z_2}=-z_1^2\frac{\partial}{\partial
z_1}+mz_2^{m-1}\zeta_2\frac{\partial}{\partial \zeta_1}
\]
and
\[
\frac{\partial}{\partial \zeta_2}=z_1^{-m}\frac{\partial}{\partial
\zeta_1},
\]
we have
\[
\frac{\partial}{\partial z_2}\wedge\frac{\partial}{\partial
\zeta_2}=-z_1^{2-m}\frac{\partial}{\partial z_1}\wedge\frac{\partial}{\partial \zeta_1}.
\]
Therefore on $(U_1\cap U_2)\times\mathbb{P}^1$,
\begin{align*}
&v_2(z_2,\zeta_2)\\
=&\Big{(}b_{2k}(z_2)\zeta_2^{2k}+b_{2k-1}(z_2)\zeta_2^{2k-1}+...+b_0(z_2)\Big{)}\Big{(}\frac{\partial}{\partial z_2}\wedge\frac{\partial}{\partial \zeta_2}\Big{)}^{\otimes k}\\
=&\bigg{[}b_{2k}(\frac{1}{z_1})z_1^{2km}\zeta_1^{2k}+b_{2k-1}(\frac{1}{z_1})z_1^{(2k-1)m}\zeta_1^{2k-1}+...+b_0(\frac{1}{z_1})\bigg{]} (-z_1^{2-m})^k\Big{(}\frac{\partial}{\partial z_1}\wedge\frac{\partial}{\partial \zeta_1}\Big{)}^{\otimes k}\\
=&(-1)^k\bigg{[}b_{2k}(\frac{1}{z_1})z_1^{(m+2)k}\zeta_1^{2k}+b_{2k-1}(\frac{1}{z_1})z_1^{(m+2)k-m}\zeta_1^{2k-1}+...+
b_0(\frac{1}{z_1})z_1^{(2-m)k}\bigg{]} \Big{(}\frac{\partial}{\partial
z_1}\wedge\frac{\partial}{\partial \zeta_1}\Big{)}^{\otimes k}
\end{align*}
Comparing the coefficients above with the coefficients of
\[
v_1(z_1,\zeta_1)=\Big{(}a_{2k}(z_1)\zeta_1^{2k}+a_{2k-1}(z_1)\zeta_1^{2k-1}+...+a_0(z_1)\Big{)}\Big{(}\frac{\partial}{\partial
z_1}\wedge\frac{\partial}{\partial \zeta_1}\Big{)}^{\otimes k},
\]
we get
\[
\left\{
  \begin{array}{l}
    (-1)^k b_{2k}(\frac{1}{z_1})z_1^{2k+km}=a_{2k}(z_1), \\
    (-1)^k b_{2k-1}(\frac{1}{z_1})z_1^{2k+(k-1)m}=a_{2k-1}(z_1),\\
    \cdots \\
    (-1)^k b_{i}(\frac{1}{z_1})z_1^{2k+(i-k)m}=a_{i}(z_1),\\
    \cdots\\
    (-1)^k b_{0}(\frac{1}{z_1})z_1^{2k-mk}=a_{0}(z_1).
  \end{array}
\right.
\]

Thus we see that for $a_i$ and $b_i$ to be holomorphic, $a_i(z_1)$ is a polynomial of degree $\le 2k+(i-k)m$ (thus it has $2k+(i-k)m+1$ degrees of freedom for $i\ge k-\frac{2k}{m}$ and $0$ otherwise), and then $b_i$ is uniquely determined by $a_i$. Therefore,
\begin{align*}
\dim H^0(-kK_{F_m})&=\dim \Gamma\big{(}\big{(}\bigwedge\nolimits^2T_{F_m}\big{)}^{\bigotimes k}\big{)}\\
&=\big{(}2k+km+1\big{)}+\big{(}2k+(k-1)m+1\big{)}+...+\big{(}2k-\Big{[}\frac{2k}{m}\Big{]}m+1\big{)}\\
&=\frac{1}{2}\Big{(}4k+(k-\Big{[}\frac{2k}{m}\Big{]})m+2\Big{)}\Big{(}k+\Big{[}\frac{2k}{m}\Big{]}+1\Big{)}.
\end{align*}
Thus Theorem \ref{Hirzebruch} is proved.

\medskip
It is then a direct calculation to get the following
\begin{cor} The dimension of the first cohomology of $mK_{F_m}$ is
$$h^1(kK_{F_m})=\frac{1}{2}\Big{(}4k-2+(k-1)m-\Big{[}\frac{2(k-1)}{m}\Big{]}m\Big{)}\Big{(}k+\Big{[}\frac{2(k-1)}{m}\Big{]}\Big{)}-4k^2+4k-1.$$
\end{cor}

\section{Deformation of a blow-up of a projective space}

Take $v$ distinct points $q_1,\cdots, q_v$ on $\mathbb{P}^n$, let
$M\to\mathbb{P}^n$ be the blow-up of $\mathbb{P}^n$ along these $v$ points, and let $E_1,\cdots,E_v$ be the corresponding exceptional
divisors.

Let $\mathbb{Z}_+$ be the set of nonnegative integers. For $\alpha=(\alpha_1,\cdots, \alpha_n)\in
\mathbb{Z}_+^n$, define $|\alpha|=\sum_{i=1}^n\alpha_i$. Let $f(z_1, \dots, z_n)$ be a polynomial
of $n$ variables. Denote
$$\frac{\partial^{|\alpha|}f}{\partial
z^\alpha}=\frac{\partial^{|\alpha|}f}{\partial
z_1^{\alpha_1}\cdots\partial z_n^{\alpha_n}}.$$

We need the following lemma:

\begin{lem} For an integer $k>0$,
 the vector space
$H^0(-kK_M)$ is isomorphic to the vector space consisting of polynomials $f(z_1,\dots, z_n)$ of
degree $\le (n+1)k$ satisfying
$$\frac{\partial^{|\alpha|}f}{\partial z^\alpha}(q_i)=0, \quad \forall 1\le
i\le v, \forall |\alpha|<(n-1)k.$$
\end{lem}
\begin{proof}
The line bundle $-kK_M$ (resp. $-kK_{\mathbb{P}^n}$) is isomorphic
to $(\bigwedge^nT_M)^{\otimes k}$ (resp.
$(\bigwedge^nT_{\mathbb{P}^n})^{\otimes k}$). Thus the natural
push-forward of vector fields $$\pi_*: \Gamma(T_M)\to
\Gamma(T_{\mathbb{P}^n})$$ induces a natural push-forward
$$\widetilde{\pi}_*: \Gamma\big{(}\big{(}\bigwedge\nolimits^nT_M\big{)}^{\otimes k}\big{)}\to
\Gamma\big{(}\big{(}\bigwedge\nolimits^nT_{\mathbb{P}^n}\big{)}^{\otimes k}\big{)}$$ hence a
natural push-forward
$$\widetilde{\pi}_*: \Gamma(-kK_M)\to
\Gamma(-kK_{\mathbb{P}^n}).$$

Since $\pi: M\to \mathbb{P}^n$ is an isomorphism outside the
exceptional divisors,  $\widetilde{\pi}_*$ must be injective.
Indeed, if a holomorphic section of $\Gamma(-kK_M)$ is zero on a
open subset of $M$, then it must be identically zero. Therefore the
dimension of $\Gamma(-kK_M)$ is equal to the image of
$\widetilde{\pi}_*$.

Now we examine which elements of $\Gamma(-kK_{\mathbb{P}^n})$ are in
the image of $\widetilde{\pi}_*: \Gamma(-kK_M)\to
\Gamma(-kK_{\mathbb{P}^n})$. Locally, we consider an open polydisc
$U=\{z=(z_1,\cdots,z_n): |z_i|<\epsilon, \forall i\}$ of
$\mathbb{C}^n$ where $\epsilon>0$ is sufficiently small. Let
$\widetilde{U}$ be the blow-up of $U$ along the center
$(0,\cdots,0)$. Then
$$\widetilde{U}=\{(z,\zeta)\in U\times \mathbb{P}^{n-1}: z_i\zeta_j=z_j\zeta_i, \forall 1\le i<j\le n\},$$
where $(\zeta_1,\cdots,\zeta_n)$ is the homogeneous coordinates of
$\zeta\in\mathbb{P}^{n-1}$. Let
$$\widetilde{U}_i=\{(z,\zeta)\in\widetilde{U}: \zeta_i\neq 0\}.$$
Without loss of generality, we consider $i=1$.
$$\widetilde{U}_1=\{(w_1,\cdots,w_n): |w_1|, |w_1w_2|, \cdots, |w_1w_n|<\epsilon\}$$
with the embedding $\widetilde{U}_1\hookrightarrow \widetilde{U}$
given by $(w_1,w_2,\cdots,w_n)\mapsto (w_1,w_1w_2,\cdots,w_1w_n)$.

An element of $(\bigwedge^nT_U)^{\otimes k}$ is of the form
$$f(z_1,z_2,\cdots,z_n)\big{(}\frac{\partial}{\partial z_1}\wedge\cdots\wedge\frac{\partial}{\partial z_n}\big{)}^{\otimes k}.$$
Since $w_1=z_1, w_2=z_2/z_1,\cdots, w_n=z_n/z_1$, we have
$\frac{\partial}{\partial z_1}=\frac{\partial}{\partial
w_1}-\frac{w_2}{w_1}\frac{\partial}{\partial
w_2}-\frac{w_3}{w_1}\frac{\partial}{\partial
w_3}-\cdots-\frac{w_n}{w_1}\frac{\partial}{\partial w_n}$,
$\frac{\partial}{\partial z_2}=\frac{1}{w_1}\frac{\partial}{\partial
w_2}$, $\cdots$,$\frac{\partial}{\partial
z_n}=\frac{1}{w_1}\frac{\partial}{\partial w_n}$. Therefore it
induces a meromorphic section of
$(\bigwedge^nT_{\widetilde{U}_1})^{\otimes k}$:
$$
f(w_1,w_1w_2,\cdots,w_1w_n)\frac{1}{w_1^{(n-1)k}}\big{(}\frac{\partial}{\partial
w_1}\wedge\cdots\wedge\frac{\partial}{\partial w_n}\big{)}^{\otimes
k}.$$ Expanding $f$ in a power series, we see that the smallest
degree of nonzero terms of the power series should be no less than
$(n-1)k$ to ensure the above section being holomorphic.

Thus, the local calculation shows  $\Gamma(-kK_M)$ is isomorphic to
the subspace of $\Gamma(-kK_{\mathbb{P}^n})$ where the Taylor
expansion of each section  at $q_i (1\le i\le v)$ contains no term
of degree $< (n-1)k$.

We can always find a hyperplane $H\subset \mathbb{P}^n$ that does not contain any $q_i (1\le i\le v)$. The complement of $H$, denoted
by $V$, is isomorphic to $\mathbb{C}^n$. It is well known that
sections of
$\Gamma(-kK_{\mathbb{P}^n})=\mathcal{O}_{\mathbb{P}^n}((n+1)k)$ are in one-to-one correspondence with polynomials on $V\cong\mathbb{C}^n$ of degree $\le (n+1)k$. Thus the dimension of
$\Gamma(-kK_{\mathbb{P}^n})$ is ${(n+1)k+n\choose n}$.  Therefore the subspace $\Gamma(-kK_M)=H^0(-kK_M)$  can be identified with the set of polynomials $f$ of degree $\le (n+1)k$ satisfying
$$
\frac{\partial^{|\alpha|}f}{\partial z^\alpha}(q_i)=0, \quad \forall 1\le
i\le v, \forall |\alpha|<(n-1)k.
$$
\end{proof}

Now we examine the case when $n=2$ and $k=1$, i.e. $H^0(-K_M)$
where $M$ is the blow-up of $\mathbb{P}^2$ at $v$ points
$q_1,\cdots, q_v$. The above lemma says that
$$
H^0(-K_M)= \{\textrm{polynomials $f$ of degree $\le 3$ such that } f(q_i)=0,
\forall 1\le i\le v.\}
$$
The dimension of $H^0(-K_M)$ depends on the position of
$q_i$'s. Indeed, let $f=\sum_{i+j\le 3}a_{ij}z_1^iz_2^j$ be the polynomial corresponding to the vector $v_f=(a_{ij})_{i+j\le 3}\in
\mathbb{C}^{10}$. Consider the map  $\mathbb{C}^2\to \mathbb{C}^{10}$ sending a point $q=(x,y)$ to $\hat{q}=(x^iy^j)_{i+j\le 3}=(1,x,y, x^2,xy,y^2, x^3,x^2y,xy^2,y^3)$ (the map $q\mapsto \hat{q}$ can be thought of as the affine version of Veronese embedding) and denote the image of this map by $S$. Then the condition $f(q)=0$ means exactly that $\hat{f}\cdot\hat{q}=0$. As there is no proper linear subspace of $\mathbb{C}^{10}$  containing $S$, generic $v$ points $\hat{q}_1,\cdots, \hat{q}_v$ in $S$ are linearly independent for $v\le 10$. Therefore, $H^0(-K_M)=10-v$ for generic $v (\le 10)$
points $p_1,\cdots, p_v$.

Now we assert that for any $v\le 4$ points $q_i=(x_i, y_i)$, the
corresponding $v$ points $\{\hat{q}_i\}_{i=1}^v$ in $S$ span a
linear subspace of dimension $v$ in $\mathbb{C}^{10}$ (in another
word, they are in general position). Indeed, after changing
coordinates if necessary, we can assume all $x_i$'s are different.
Then the matrix with the $i$-th row being the coordinate of
$\hat{q}_i$ is
$$\left(
                                               \begin{array}{cccccccccc}
                                                 1 & x_1 & y_1 & x_1^2 & x_1y_1 & y_1^2 & x_1^3 & x_1^2y_1 & x_1y_1^2 & y_1^3 \\
                                                 \vdots &  \vdots &\vdots &\vdots &\vdots &\vdots &\vdots &\vdots &\vdots &\vdots \\
                                                 1 & x_v & y_v & x_v^2 & x_vy_v & y_v^2 & x_v^3 & x_v^2y_v & x_vy_v^2 & y_v^3 \\
                                               \end{array}
                                             \right),$$
which contains a $v\times v$ square $$\left(
                                      \begin{array}{cccc}
                                        1 & x_1 & \cdots & x_1^{v-1}\\
                                        \vdots & \vdots & &\vdots \\
                                        1 & x_v & \cdots & x_v^{v-1}\\
                                      \end{array}
                                    \right)$$
Its determinant is $\prod_{1\le i<j\le v}(x_j-x_i)$ which is nonzero
since we assume that the $x_i$'s are different. Therefore the above $v\times 10$ matrix is of full rank, which implies the $v$ points
$\{\hat{q}_i\}_{i=1}^v$ in $S$ span a linear subspace of dimension
$v$ in $\mathbb{C}^{10}$.

On the other hand, we can easily find as finitely many points as we want in $S$ which lie in a $4$
dimensional subspace of $\mathbb{C}^{10}$. For example, we can take $\{q_i=(i,0)\}_{i=1}^v$, the
corresponding $\{\hat{q}_i\}_{i=1}^v\in\mathbb{C}^{10}$ satisfy the requirement. for $v\ge 5$,
start from $v$ points of $S$ in a $4$ dimensional subspace, move them point by point to a general
position, we get $v$ points in $S$ spanning a subspace of dimension $4, 5, \cdots, \min(v,10)$,
respectively. Therefore we have proved the following theorem:

\begin{thm}
Let $M\to\mathbb{P}^2$ be the blow-up of $\mathbb{P}^2$ along $v$ distinct points. The dimension of $H^0(-K_M)$ is $10-v$ for $v\le
4$; for $v>4$, the dimension depends on the position of the $v$
points: for each integer $i$ such that $\max(10-v,0)\le i\le6$,
there exist $v$ distinct points such that the corresponding $M$
satisfies $\dim H^0(-K_M)=i$.
\end{thm}

\medskip
It is then a direct calculation to get the following
\begin{cor} The dimension of the first cohomology of $2K$ is
$$h^1(2K)=\left\{
    \begin{array}{ll}
      0, & \hbox{if $v\le 4$;} \\
      \hbox{integer between $max(0,v-10)$ and $v-4$}, & \hbox{if $v\ge 5$.}
    \end{array}
  \right.$$

  \end{cor}

\bigskip

\noindent {\sc Ning Hao }\\
Mathematics Department, SUNY at Stony Brook\\
Stony Brook, NY 11794, US\\
Email: {\tt nhao@math.sunysb.edu}
\bigskip

\noindent {\sc Li Li }\\
Korean Institute for Advanced Study\\
207-43 Cheongryangri-dong, Seoul 130-722, Korea\\
Email: {\tt lili@kias.re.kr}\\

\end{document}